\documentclass[10pt]{amsart}

\newtheorem{theorem}{Theorem}[section]
\newtheorem{proposition}{Proposition}[section]
\newtheorem{lemma}{Lemma}[section]
\newtheorem{corollary}{Corollary}[section]

\theoremstyle{definition}
\newtheorem{remark}{Remark}[section]
\newtheorem{question}{Question}[section]

\numberwithin{equation}{section}

\newcommand{\Nat}{{\mathbb N}}

\newcommand{\Real}{{\mathbb R}}
\newcommand{\Com}{{\mathbb C}}

\newcommand{\A}{{\mathcal A}}
\newcommand{\B}{{\mathcal B}}
\newcommand{\C}{{\mathcal C}}

\begin{document}

\title[Arithmetical Aspects of Beurling's Theorem]
{Arithmetical Aspects of Beurling's Real Variable Reformulation of the Riemann
Hypothesis}

\author{Luis B\'{a}ez-Duarte}

\date{20 October 2000}

\email{lbaez@ccs.internet.ve}

%\thanks{This paper is not for wide distribution or publication in its present form}

\begin{abstract}
Let $\rho(x):=x-[x]$, $\chi:=\chi_{(0,1]}$, the characteristic function of $(0,1]$, $\lambda(x):=\chi(x)\log x$,
and 
$M(x):=\sum_{k\leq x}\mu(k)$, where $\mu$ is the M\"{o}bius function. $\B$ is the space of functions defined in
$(0,\infty)$ by expressions
$\sum_{k=1}^{n}c_k\hspace{.3mm}\rho(\theta_k/x)$ with $n\in\Nat$,
$c_k\in\Com$ and $\theta_k\in(0,1]$. A minor sharpening of the results of B. Nyman and A. Beurling states that for
any fixed $p\in(1,\infty)\overline{\B}^{L_p}$ the Riemann zeta function $\zeta(s)\not=0$ for $\Re s>1/p$, if and
only if $L_p(0,1)\subset \overline{\B}^{L_p}$, which, furthermore, is equivalent to
$\chi\in\overline{\B}^{L_p}$ or to $\lambda\in\overline{\B}^{L_p}$. Starting from the elementary identity
$\lambda(x):=\int_{0}^1M_1(\theta)\rho(\theta/x)\theta^{-1}d\theta$, with $M_1(\theta):=M(1/\theta)$, where
the integral suggests a limit of functions in $\B$, we were led to the following two
\textit{arithmetical} versions of the Nyman-Beurling results, proved by classical, quasi elementary,
number-theoretic methods. Define $G_n$, a \textit{natural approximation} to $\lambda$, by
$G_n(x):=\int_{1/n}^1M_1(\theta)\rho(\theta/x)\theta^{-1}d\theta$, then for all $p\in(1,\infty)$\\
\ \\
(I) $\|G_n-\lambda\|_p\rightarrow0$ implies $\zeta(s)\not=0$ in $\Re s\geq1/p$, and $\zeta(s)\not=0$ in 
$\Re s>1/p$ implies $\|G_n-\lambda\|_r\rightarrow0$ for all $r\in(1,p)$.\\
\ \\
Likewise noting that $\zeta(s)\not=0$ in $\Re
s>1/p$ is equivalent to $\|M_1\|_r<\infty$ for all $r\in(1,p)$, we have for all $p\in(1,\infty)$\\
\ \\  
(II) $\|M_1\|_p<\infty$ implies $\lambda\in\overline{\B}^{L_p}$, and $\lambda\in\overline{\B}^{L_p}$ implies
$\|M_1\|_r<\infty$ for all $r\in(1,p)$.\\
\ \\
It is clear from (I) that $G_n\rightarrow\lambda$ diverges in $L_2$, although it is shown to converge both pointwise and in
$L_1$ to $\lambda$. The general $L_p$ case is also discussed. Some older \textit{natural
approximations} to $\chi$, for which J. Lee, M. Balazard and E. Saias proved theorems
analogous to (I), are shown to diverge in $L_2$.

\end{abstract}

\maketitle

\section{Introduction}
\subsection{Preliminaries and notation}
For every $p\in[1,\infty]$ we canonically imbed $L_p(0,1)$ in $L_p(0,\infty)$. The conjugate
index is always denoted by $q:=p/(p-1)$. $\rho(x):=x-[x]$ stands throughout for the fractional part of
the real number $x$, and $\chi:=\chi_{(0,1]}$ is the characteristic function of
the set $(0,1]$. We define the function $\lambda$ by
\begin{equation}\label{lambda}
\lambda(x):=\chi(x)\log x.
\end{equation}
For every $a>0$ the operator $K_a$ given by $$K_a f(x):=f(ax),$$ acts continuously on every
$L_p(0,\infty)$ to itself, for $1\leq p \leq \infty$. 

$\A$ shall be the vector space of functions $f$ of the form
\begin{equation}\label{abcfunction}
f(x)=\sum_{k=1}^{n}c_k \rho\left(\frac{\theta_k}{x}\right),
\end{equation}
with $n\in\Nat$, $c_k\in\Com$, $\theta_k > 0$, $1\leq k \leq n$. For $E\subseteq (0,\infty)$
denote by $\A_E$ the subspace of $\A$ where the $\theta_k\in E$. In particular we let
$\B=\A_{(0,1]}$. $\C$ is the subspace of $\B$ resulting from requiring that
\begin{equation}\label{sumckthetakiszero}
\sum_{k=1}^{n}c_k\theta_k=0.
\end{equation}
Clearly $$\C\subset\B\subset\A\subset L_p(0,\infty)$$ for $1<p\leq\infty$. Note that functions
in $\C$ vanish in $(1,\infty)$, so $\C\subset L_p(0,1)$ for $1\leq p\leq\infty$. $\A$ is invariant under
any $K_a$, $a>0$, while $\B$ and $\C$ are invariant under $K_a$ for $a\geq 1$.

Recall the usual aritmetical functions $M$ and $g$ given by
$$M(x)=\sum_{k\leq x}\mu(k),$$
\begin{equation}\label{g}
g(x):=\sum_{k\leq x}\frac{\mu(k)}{k},
\end{equation}
where $\mu$ is the arithmetical M\"{o}bius function. We shall denote
\begin{equation}\label{M1}
M_1(\theta):=M\left(\frac{1}{\theta}\right).
\end{equation}
It is classical number theory that both $M(x)x^{-1}\rightarrow0$ and $g(x)\rightarrow0$ as
$x\rightarrow\infty$ are \textit{elementarily}\footnote{Heretofore \textit{elementary} is to be
understood in the traditional number theoretical sense, as ``no analytic function theory", ``no
Fourier analysis".} equivalent to the prime number theorem. A stronger but still elementary estimate
is $M(x)\ll x(\log x)^{-2}$.
 
Let us also define the less common
$\gamma$ and
$H_p$ by
\begin{equation}\label{gamma}
\gamma(x):=\sum_{k\leq x-1} \frac{M(k)}{k(k+1)},
\end{equation}
\begin{equation}\label{Hp}
H_p(x)=\int_1^x M(t)t^{-2/p}dt.
\end{equation}
Summing (\ref{g}) by parts we get
\begin{equation}\label{Mggammanatural}
g(n)=\frac{M(n)}{n}+\gamma(n),\ \ (n\in\Nat),
\end{equation}
and trivially from $|M(x)|\leq x$
\begin{equation}\label{Mggammareal}
g(x)=\frac{M(x)}{x}+\gamma(x) +O(1/x),\ \ (x\in\Real).
\end{equation}

\subsection{The weak Nyman-Beurling theorem}
An easy consequence of Wiener's $L_2$ Tauberian theorem (cfr. \cite{rudin}) is that $\A$ is dense in
$L_2(0,\infty)$ (see \cite{IUO}). B. Nyman \cite{nyman} for $L_2$ and A.
Beurling \cite{beurling} for general $L_p$ obtained the much deeper result:

\begin{theorem}[Nyman-Beurling]\label{nymanbeurlingtheorem}
The Riemann zeta-function is free from zeroes in the half-plane $\sigma>1/p$, $1<p<\infty$, if and
only if $\C$ is dense in the space $L_p(0,1)$, which is equivalent to $-\chi\in\overline{\C}^{L_p}$.
\end{theorem}
\noindent To prove this theorem Beurling first noted that $\C$ is dense in $L_p$ if and only
$-\chi\in\overline{\C}^{L_p}$, then showed quite simply that $-\chi\in\overline{\C}^{L_p}$ implies
$\zeta(s)\not=0$ for $\Re s>1/p$. The proof of the converse, which, in his own words, is \textit{less trivial}, is
by contradiction. If $-\chi\not\in\overline{\C}^{L_p}$, then, of course, $\C$ is not dense in the space
$L_p(0,1)$. But this, by a highly involved functional analysis argument, implies the existence of a zero
with real part greater than
$1/p$. Later proofs of this fact are illuminating, but just as difficult (see \cite{donoghue}, \cite{bercovici},
\cite{vienna}). The degree to which the apparent depth of the two sides of the proof is so starkly contrasting
has led some authors to voice doubts about the usefulness of the Nyman-Beurling approach (see, for example,
\cite{levinson}), yet, it has led others to attempt to level off the two sides of the proof.\\  
\ \\
We say that $\phi$ is a \textit{generator}\footnote{We called these generators \textit{ strong generators} in \cite{IUO},
and applied the term \textit{generator} when $a$ was allowed to range in $(0,\infty)$ in (\ref{definegenerator}).} if
$\phi\in L_p(0,\infty)$ for all
$p\in(1,\infty)$ and 
\begin{equation}\label{definegenerator}
L_p(0,1)\subseteq\text{span}_{L_p}\{K_a \phi\}_{a\geq1}, \ \ \ (1<p<\infty).
\end{equation}
$-\chi$ is the simplest example of a generator (the minus sign is immaterial, but more convenient). The
function
$\lambda$ defined in (\ref{lambda}) is also a generator since
\begin{equation}\label{lambdaimplieschi}
\frac{1}{a-1}(K_a - I)\lambda \overset{L_p}{\rightarrow} \chi,\ \ \ (a\downarrow1),\ \ \ (1\leq p
<\infty).
\end{equation}
Clearly any generator $\phi$ may well take the place of $-\chi$ in Theorem \ref{nymanbeurlingtheorem}.
These considerations, together with the fact that
$$
f(x)=\frac{1}{x}\sum_{k=1}^{n}c_k\theta_k,\ \ \ (x>1),
$$
for every $f\in\B$ as in (\ref{abcfunction}), allow the following minor extension of the Nyman-Beurling
Theorem \ref{nymanbeurlingtheorem}, where reference to density of $\C$ or $\B$ is dropped. 

\begin{theorem}\label{extendednymanbeurlingtheorem}
Let $\phi$ be a generator and $p\in(1,\infty)$. Then $\zeta(s)\not=0$ for $\Re s > 1/p$ if and
only if $\phi\in\overline{\B}^{L_p}$.
\end{theorem}
\noindent Obviously the above theorem implies this weaker version:

\begin{theorem}[Weak Nyman-Beurling Theorem]\label{weaknymanbeurlingtheorem}
Let $p\in(1,\infty)$ and $\phi$ be a generator. Then $\zeta(s)\not=0$ for $\Re s
> 1/p$ if and only if $\phi\in\overline{\B}^{L_r}$ for $r\in(1,p)$.
\end{theorem}

\noindent Direct, independent proofs of this theorem for $\phi=-\chi$, not depending on deep functional
analysis results were achieved independently by J. Lee \cite{lee}, and M. Balazard and E. Saias
\cite{notes1}. These proofs only make use of standard number theoretical techniques. Thus Lee,
not inapropriately, presents his result \textit{an arithmetical version} of Beurling's theorem. The
\textit{only if} part of these proofs depends on identifying \textit{natural approximation}s $f_n$,
which we define as sequences in $\C$ or $\B$, such that this \textit{weak implication} holds for all $p\in(1,\infty)$: 
\begin{equation}\label{weakimplication}
(\zeta(s)\not=0,\  (\Re s > 1/p)) \Rightarrow (\|f_n-\phi\|_r\rightarrow0,\ \forall r\in(1,p)).
\end{equation}
Balazard and Saias \cite{notes1} asked the natural question:

\begin{question}\label{weaktostrong} 
For a given specific natural approximation $\{f_n\}$ is it true for some or all $p\in(1,\infty)$ that
the weak implication (\ref{weakimplication}) can be substituted for the \textit{strong implication}
\begin{equation}\label{strongimplication}
(\zeta(s)\not=0,\  (\Re s > 1/p)) \Rightarrow (\|f_n-\phi\|_p\rightarrow0)?
\end{equation}
\end{question}
\noindent We shall answer this question mostly in the negative in Section \ref{divergenceofnaturalapproximations}. 
The first such natural approximation $\{B_n\}\subset\B$ had appeared earlier in
\cite{BR1} defined by

\begin{equation}\label{Bn}
B_n(x):=\sum_{k=1}^{n}\mu(k)\rho\left(\frac{1}{kx}\right)-ng(n)\rho\left(\frac{1}{nx}\right).
\end{equation}
This sequence arises rather naturally in more than one way: it is the unique answer to the problem
of finding
$f\in\C$ as in (\ref{abcfunction}) with $\theta_k=1/k$, and $f(k/n)=-1$ for $1\leq k\leq n-1$. Or it can also
be seen as a truncation of the fundamental identity
\begin{equation}\label{natseries}
-1=\sum_{k=1}^{\infty}\mu(k)\rho\left(\frac{1}{kx}\right), \ \ \ (x>0),
\end{equation}
It is easily seen that
$B_n(x)=-1$ in $[1/n,1]$, and using the prime number theorem we proved that
\begin{equation}\label{BnconvergesinL1}
\|\chi+B_n\|_1\rightarrow 0,
\end{equation}
which led us to ask whether the strong or the weak implications (\ref{weakimplication}),
(\ref{strongimplication}) were true for $f_n=B_n$, $1<p\leq 2$. A mild positive answer was (\cite{BR1}, Proposition
2.4) that $\zeta(s)$ has a non-trivial zero-free half-plane if and only if
$\|\chi+B_n\|_p\rightarrow 0$ for some $p>1$, which conferred some legitimacy to the question. In
related work V. I. Vasyunin \cite{vasyunin}, referring to earlier results of N. Nikolski
\cite{nikolski}, took up the study of the $L_2$ case in quite some depth for a $B_n$-related sequence
$\{V_n\}\subset\C$ defined by
\begin{equation}\label{Vn}
V_n(x):=\sum_{k=1}^{n}\mu(k)\rho\left(\frac{1}{kx}\right)-g(n)\rho\left(\frac{1}{x}\right).
\end{equation}
Vasyunin also conducted numerical studies leading him to state that \textit{we can hardly hope that
the series converges in the $L_2$-norm}. That this is indeed the case was first proved in \cite{IUO}. The sequence
$\{S_n\}\subset\B$ defined by
\begin{equation}\label{Sn}
S_n(x):=\sum_{k=1}^{n}\mu(k)\rho\left(\frac{1}{kx}\right),
\end{equation}
perhaps the most \textit{natural} in view of (\ref{natseries}), is obviously $L_2$-equivalent to $\{V_n\}$ since
$g(n)\rightarrow 0$. The relationship with $B_n$ is more complicated, however, since by Corollary
\ref{arithmeticalfunctionbehavior} below the $L_p$-norm of $ng(n)\rho(1/nx)$ is of order $|g(n)|n^{1/q}$ which 
does not tend to zero if $\zeta(s)$ has a zero with real part $1/p$, such being the case, of course if $p=2$. Furthermore
$B_n$ is not a \textit{series} as defined in (\ref{seriesinC}) while $V_n$ is the most natural series. 

J. Lee \cite{lee} proved the weak Theorem \ref{weaknymanbeurlingtheorem} using $V_n$, $1<p\leq2$, and, independently,
M. Balazard and E. Saias \cite{notes1} did likewise for $B_n$ and $S_n$, $1<p<\infty$.

A further approximating sequence $\{F_n\}\subset\C$ promoted in \cite{BR1} as the \textit{dual
approximation}, given by
\begin{equation}\label{Fn}
F_n(x):=
\sum_{k=1}^n
\left(M\left(\frac{n}{k}\right)-M\left(\frac{n}{k+1}\right)\right)\rho\left(\frac{k}{nx}\right)
-\rho\left(\frac{1}{nx}\right),
\end{equation}
is of a different nature, as the $\theta_k$ are uniformly distributed in $(0,1)$ as
$n\rightarrow\infty$. It is proved in \cite{BR1} that $\|F_n+\chi\|_1\rightarrow0$, and it can also
be shown that $F_n(x)\rightarrow-1$ for $0<x\leq1$. The following question is however open:

\begin{question}
Is $F_n$ a natural approximation? 
\end{question}

\subsection{Description of main results.}

The purpose of this paper is twofold. In first place we produce in Section \ref{arithmeticaltheorems} two
\textit{ arithmetical} versions of the Nyman-Beurling results, standing somewhere in between the
strong and the weak Theorems \ref{extendednymanbeurlingtheorem} and \ref{weaknymanbeurlingtheorem}. We state them
as theorems A and B below, and prove them later as Theorems \ref{maintheoremversion1},
\ref{maintheoremversion2}. To discuss them properly we mention first the following proposition, a version of
Littlewood's criterion for the Riemann hypothesis, established below as Proposition
\ref{arithmeticalRHp}:\\
\ \\
\noindent \textbf{Proposition.}\textit{ For all $p\in(1,\infty)$, $\zeta(s)\not=0$ for $\Re s > 1/p$
if and only if
$\|M_1\|_r<\infty$ for all $r\in(1,p)$,}\\
\ \\
\noindent and introduce a new ``natural approximation", of more general type, $G_n\in\overline{\A_{(1/n,1]}}^{L_p}$,
$p\in(1,\infty)$, defined by
\begin{equation}\label{Gn}
G_n(x):=\int_{1/n}^{\infty}M_1(\theta)\rho\left(\frac{\theta}{x}\right)\frac{d\theta}{\theta},
\end{equation} 
arising, among others in Section \ref{divergenceofnaturalapproximations}, from the convolution
(\ref{mobiusidentity2})
$$\lambda(x)=\int_0^1 M_1(\theta)\rho\left(\frac{\theta}{x}\right)\frac{d\theta}{\theta}.$$

\noindent Our two main theorems are then:\\
\ \\
\noindent \textbf{Theorem A.} (Arithmetical Nyman-Beurling Theorem, I). \textit{The following
statements are true \textit{for all} $p\in(1,\infty)$.}\\
\ \\
(a) $\|M_1\|_p<\infty$ \textit{implies} $\lambda\in\overline{\B}^{L_p}$.\\
\ \\
(b) $\lambda\in\overline{\B}^{L_p}$ \textit{implies} $\|M_1\|_r<\infty$ \textit{for all} $r\in(1,p)$.\\
\ \\
\noindent \textbf{Theorem B.} (\textit{Arithmetical Nyman-Beurling Theorem, II\textnormal{)}}
\textit{The following statements are true \textit{for all} $p\in(1,\infty)$}.\\
\ \\
(c) $\zeta(s)\not=0$, $\Re s>1/p$ implies $\|G_n-\lambda\|_r\rightarrow0$ for all $r\in(1,p)$.\\
\ \\
(d) $\|G_n-\lambda\|_p\rightarrow0$ implies $\zeta(s)\not=0$, $\Re s\geq 1/p$.\\
\ \\
\noindent We give rather simple proofs of these statements. Actually, elementary for (a), and
quasi elementary for (c). Further note that (a) and (d) are \textit{strong} statements.\\
\ \\ 
Secondly, in section \ref{divergenceofnaturalapproximations}, the paper aims to explore the delicate gap
between the weak and strong forms of the Nyman-Beurling theorem. We shall show that all the \textit{natural
approximations} $B_n$, $V_n$, $S_n$, $F_n$, and $G_n$ diverge in $L_2$. We also study the general $L_p$ case. The most
interesting conclusion is this: If the Riemann hypothesis were not true, and $1/p=\sup\{\Re s|\hspace{1mm}\zeta(s)=0\}$,
then $S_n$, $V_n$, and $G_n$ would also diverge in $L_p$ provided there is a zero of real part $1/p$. We have not decided
the question for $B_n$ and $F_n$.

\section{Technical lemmae and preliminary propositions}\label{lemmae}
Throughout this section $1<p<\infty$. Some of the results herein may be part of the common
folklore and/or stated in less general form than is possible. They are listed here however for the
sake of completeness and readability. We thank A. M. Odlyzko for his generous help in these matters.\\
\subsection{Technical Lemmae}
It is assumed that $f$ is a locally bounded complex
valued function defined on $[1,\infty)$, whose Mellin transform  $\tilde{f}$, defined here as
\begin{equation}\label{defmellin}
\tilde{f}(s):=\int_1^{\infty}f(x)x^{-s-1}dx,
\end{equation} has a finite abcissa of convergence $\alpha=\alpha_{f}$.  

\begin{lemma}[Order Lemma]\label{orderlemma}
If $\tilde{f}(s)$ has a pole at $s_0=\sigma_0+it_0$ in a meromorphic extension to a possibly larger
half-plane, then $f(x)\not=o\left(x^{\sigma_0}\right)$.
\end{lemma}
\begin{proof}
This is just an adaptation of the proof of $M(x)\not=o(\sqrt x)$ in
\cite{titchmarsh}. Assume by contradiction that
$f(x)=o\left(x^{\sigma_0}\right)$. A fortiori
$f(x)=O\left(x^{\sigma_0}\right)$, so that the integral in (\ref{defmellin}) would
actually converge in $\Re s>\sigma_0$. Now let $s=\sigma+it_0$ with
$\sigma\downarrow\sigma_0$. If $m\geq1$ is the order of the pole, then we have
\begin{equation}\label{asymptotic}
\tilde{f}(\sigma+it_0)\sim\frac{C}{(\sigma-\sigma_0)^m}, \ \ (\sigma\rightarrow\sigma_0)
\end{equation}
for some $C\not=0$. On the other hand the little $o$ condition implies there is an $A>1$ such
that $|f(x)|<(|C|/2)x^{\sigma_0}$ for $x>A$, so that splitting the right-hand side
integral in (\ref{defmellin}) as $\int_1^A +\int_A^{\infty}$ we obtain 
\begin{equation}\nonumber
|\tilde{f}(\sigma+it_0)|\leq O_A(1) + \frac{|C|}{2(\sigma-\sigma_0)},
\end{equation}
which contradicts (\ref{asymptotic}).
\end{proof}

\begin{lemma}[Oscillation Lemma]\label{oscillationlemma}
Let $f$ be real valued. If $\alpha=\alpha_f$ is not a singularity of
$\tilde{f}(s)$ then for any
$\epsilon>0$
\begin{eqnarray}
\limsup_{x\rightarrow\infty}f(x)x^{-\alpha+\epsilon}&=&+\infty,\\\label{first}
\liminf_{x\rightarrow\infty}f(x)x^{-\alpha+\epsilon}&=&-\infty.
\end{eqnarray}
In particular, $f(x)$ changes sign an infinite number of times as $x\rightarrow\infty$.
\end{lemma}
\begin{proof}
It is obviously enough to deal with only one of the above relations. So assume that
(\ref{first}) is false. Then for some $\epsilon>0$ there is a $C$ such that
\begin{equation}\nonumber
Cx^{\alpha-\epsilon}-f(x)\geq 0. \ \ \ (x>1).
\end{equation}
Therefore
\begin{equation}\nonumber
\int_{1}^{\infty}\left(Cx^{\alpha-\epsilon}-f(x)\right)x^{-s-1}dx
=\frac{B}{s-\alpha+\epsilon}-\tilde{f}(s)
\end{equation}
is not singular at $s=\alpha$, but clearly $\alpha$ is also the abcissa of convergence of the
left-hand side integral above, which contradicts the theorem that the Laplace transform of a
positive measure has a singularity on the real axis at the abcissa of convergence
(\cite{widder}, theorem 5.b).
\end{proof}

\begin{lemma}\label{divergentintegrallemma}
Let $F:[1,\infty)\rightarrow\Com$ be locally integrable. If $\int_1^x F(t)dt \not=o(x^{1/q})$,
then $\|F\|_p=\infty$. 
\end{lemma}
\begin{proof}
It is obviously enough to consider that $F\geq 0$. By hypothesis there exists some
$\epsilon>0$, and an unbounded set $E\subset[1,\infty)$ such that 
\begin{equation}\nonumber
\int_1^y F(t)dt>\epsilon y^{1/q}, \ \ \ (\forall y\in E).
\end{equation}
Now take an arbitrary $x>1$. It is easy to see that there exists $y\in E$ such that $y>x$ and
\begin{equation}\nonumber
2\int_1^x F(t)dt<\epsilon y^{1/q}<\int_1^y F(t)dt,
\end{equation}
so that
\begin{equation}\label{step}
\int_x^y F(t)dt>\frac{\epsilon}{2} y^{1/q}.
\end{equation}
But H\"{o}lder's inequality gives
\begin{equation}\nonumber
\int_x^y F(t)dt
\leq (y-x)^{1/q}\left(\int_x^y \left(F(t)\right)^p dt\right)^{1/p},
\end{equation}
which, introduced in (\ref{step}), yields
\begin{equation}\nonumber
\int_x^y (F(t))^p dt
>\frac{\epsilon}{2}\left(\frac{y}{y-x}\right)^{1/q}>\frac{\epsilon}{2}.
\end{equation}
\end{proof}

\subsection{Some preliminary propositions}

\begin{proposition}\label{mellintransforms}
The following Mellin transforms are valid at least in the half-planes indicated.
\begin{eqnarray}\nonumber\\\label{mellinofM}
\int_1^\infty M(x)x^{-s-1}dx&=&\frac{1}{s\zeta(s)},\ \ (\Re s>1)\\\label{mellinofg}
\int_1^{\infty}(xg(x))x^{-s-1}dx&=&\frac{1}{(s-1)\zeta(s)},\ \ (\Re s>1)\\\label{mellinofgamma}
\int_{1}^{\infty}(x\gamma(x))x^{-s-1}dx&=&\frac{1}{s(s-1)\zeta(s)}+\omega(s),\ \ (\Re
s>1)\\\label{mellinofHp}
\int_1^{\infty} H_p(x)x^{-s-1}dx&=&\frac{1}{s(s+2/p-1)\zeta(s+2/p-1)},\ \ (\Re s >2/q),
\end{eqnarray}
where $\omega(s)$ is analytic in $\Re s>0$.
\end{proposition}
\begin{proof}
As in Titchmarsh's monograph \cite{titchmarsh} we write for $\Re s>1$
\begin{eqnarray}\nonumber
\frac{1}{\zeta(s)}&=&\sum_{n=1}^{\infty}(M(n)-M(n-1))n^{-s}\\\nonumber
  &=&\sum_{n=1}^{\infty}M(n)\left(n^{-s}-(n+1)^{-s}\right)\\\nonumber
  &=&s\sum_{n=1}^{\infty}M(n)\int_n^{n+1}x^{-s-1}dx\\\nonumber
  &=&s\int_1^{\infty}M(x)x^{-s-1}dx. 
\end{eqnarray}
This proves (\ref{mellinofM}). Proceed likewise with
$$\frac{1}{\zeta(s)}=\sum_{n=1}^{\infty}(g(n)-g(n-1))n^{-s+1},\ \ (\Re s>1),$$
to obtain (\ref{mellinofg}). Now using the relation (\ref{Mggammareal}) between $M(x)$, $g(x)$ and
$\gamma(x)$ subtract the preceding two Mellin transforms to get (\ref{mellinofgamma}). Finally, from
the definition (\ref{Hp}) and the trivial $|M(x)|\leq x$ we deduce $H_p(x)\ll x^{2/q}$. Next note that
$H_p$ is continuous and piecewise differentiable, which justifies the following integration by parts
at least for $\Re s > 2/q$
$$\int_1^{\infty}H_p(x)x^{-s-1}dx=\frac{1}{x}\int_1^{\infty}M(x)x^{-s-2/p}dx.$$
Now apply (\ref{mellinofM}) to arrive at (\ref{mellinofHp}).
\end{proof}
An immediate consequences of the above Mellin transforms, and the order and oscillation lemmae \ref{orderlemma},
\ref{oscillationlemma}:
\begin{corollary}\label{arithmeticalfunctionbehavior}
Each one of the functions $M(x), g(x), \gamma(x), H_p(x)$, change sign infinitely often as
$x\rightarrow\infty$. Furthermore, if $\zeta(s)$ has some zero on the line $\Re s=1/p$ then $M(x)\not=o(x^{1/p})$,
$g(x)\not=o(x^{-1/q})$, $\gamma(x)\not=o(x^{-1/q})$, and $H_p(x)\not=o(x^{1/q})$.   
\end{corollary}
\begin{remark}
In particular, $M(x)\not=o(\sqrt{x})$ (see \cite{titchmarsh}). Sharper results are of course known, e.g., that the
Mertens hypothesis is false, with $M(x)$ oscillating beyond $\pm\sqrt{x}$. This was proven by A. M. Odlyzko and
H. te Riele \cite{odlyzko}.
\end{remark}
Some further properties of $\gamma(x)$ needed later are gathered here:
\begin{lemma}\label{furtherongamma}
The function $\gamma$ satisfies\\
\ \\
(i) $\gamma(n)\rightarrow0$ as $n\rightarrow\infty$,\\
\ \\
(ii) $\gamma(n)=\int_1^n M(t)t^{-2}dt$, ($n\in\Nat$),\\
\ \\
(iii) $\int_1^{\infty} M_1(\theta)d\theta=\int_1^{\infty} M(t)t^{-2}dt=0$. This integral converges absolutely.
\end{lemma}
\begin{remark}
In \cite{PNTequivlbd} we showed that the existence of $\lim_{n\rightarrow\infty}\gamma(n)$ is elementarily
equivalent to the prime number theorem.
\end{remark}
\begin{proof}[Proof of lemma \ref{furtherongamma}]
The prime number theorem and (\ref{Mggammanatural}) yield (i). Decomposing the integral in (ii) in the intervals $(k,k+1)$
one gets (ii). Letting $n\rightarrow\infty$ in (ii) yields (iii). The absolute convergence follows from the
elementary estimate $M(x)\ll x(\log x)^{-2}$.
\end{proof}

The result on $H_p(x)$ in Corollary \ref{arithmeticalfunctionbehavior} begets some important consequences for
the norms of $M_1$.

\begin{proposition}\label{normpMinfinite}
If $\zeta(s)$ has a zero on the line $\Re s=1/p$ , then
\begin{equation}
\|M_1\|_p= \infty.
\end{equation}
\end{proposition}
\begin{remark}
Note therefore that
\begin{equation}\label{normpMfinite}
\|M_1\|_p< \infty\ \textnormal{implies (}\zeta(s)\not=0\ \textnormal{for } \Re s\geq 1/p\textnormal{)}.
\end{equation}
\end{remark}
\begin{proof}[Proof of Proposition \ref{normpMinfinite}]
Take $F(x):=M(x)x^{-2/p}$. Then $F(x)\not=o(x^{1/q})$ by Corollary \ref{arithmeticalfunctionbehavior}, so the
divergent integral lemma \ref{divergentintegrallemma} yields
$$\|M_1\|_p^p=\int_1^{\infty}|M(x)|^p x^{-2} dx = \infty.$$
\end{proof}

\begin{remark}
Since $\zeta(s)$ has roots in the critical line the above corollary tells us that
\begin{equation}\label{I2isinfinity}
\|M_1\|_2 = \infty.
\end{equation}
Using far more refined techniques S. V. Konyagin and A. Yu. Popov
\cite{konyagin} have shown a stronger result in the case $p=2$, namely
$$\int_1^x \left|M(t)\right|^2 t^{-2}dt \gg \log x. $$ 
\end{remark}

\begin{proposition}\label{arithmeticalRHp}
For any $p\in(1,\infty)$ the following statements are equivalent.\\
\ \\
(i) $\zeta(s)\not=0$ for $\Re s > 1/p$,\\
\ \\
(iv) $\|M_1\|_p <\infty$ for all $r\in(1,p).$
\end{proposition}

\begin{proof}
An extension of Littlewood's well-known criterion for the Riemann hypothesis is that condition (i) is
equivalent to 
$$M(x)\ll x^{1/p'} \textnormal{\ for all\ }p'\in(1,p),$$
(see \cite{blanchard}, proposition IV.21), so choose $p'$ with $r<p'<p$ and it is obvious how (i)
implies (iv). Now we prove that not (i) implies not (iv). So assume there is an $s_0$ with $\Re s_0
=1/p_1>1/p$ and $\zeta(s_0)=0$. Then by Corollary \ref{normpMinfinite} 
$$\|M_1\|_{p_1} =\infty,$$
but (probability space) $\|M_1\|_{p_1}\leq\|M_1\|_r$ since $p_1<r<p$.    
\end{proof}
Define the \textit{Riemann abcissa} $\beta$ by
\begin{equation}\label{RHbeta}
\beta:=\sup_{\zeta(s)=0}\Re s.
\end{equation}
Nothing is known beyond $1/2\leq\beta\leq1$. We do know however that, on the one
hand there are no zeroes on the line $\Re s=1$ and $\|M_1\|_1<\infty$, since
$M(t)\ll t(\log t)^{-2}$, and, on the other hand there are zeroes on the line $\Re s=1/2$ and
$\|M_1\|_2=\infty$. One could rightly ask the question:

\begin{question}\label{betaconjecture}
For $\beta\in(\frac{1}{2},1)$, is it true that $\|M_1\|_{1/\beta}<\infty$ if and only $\zeta(s)\not=0$ for $\Re s=\beta$. 
\end{question}

\section{Two arithmetical versions of the Nyman-Beurling Theorem}\label{arithmeticaltheorems}
 We define an operator $T$ acting on all $L_p(0,\infty)$, $p\in(1,\infty)$, by
\begin{equation}\label{Toperator}
Tf(x):=\int_0^{\infty}f(\theta)\rho\left(\frac{\theta}{x}\right)\frac{d\theta}{\theta},
\end{equation}
noting that the above integral converges
absolutely for $f\in L_p(0,\infty)$ by H\"{o}lder's inequality. Now we show that $T$ is of type $(p,p)$. This does
not follow, as could be expected, from the convolution form of the operator, on account of the difference between
the measures $dx$ and $dx/x$ in $(0,\infty)$. 

\begin{lemma}\label{Tispp}
For every $p\in(1,\infty)$ the operator $T$ is a continuous operator from $L_p(0,\infty)$ to itself.
\end{lemma}
\begin{proof}
Let $f\geq0$ and $x>0$, then splitting the range of integration at $x$ in (\ref{Toperator}) we get
\begin{equation}\nonumber
Tf(x)\leq
\frac{1}{x}\int_0^x f(\theta)d\theta +\int_x^{\infty}f(\theta)\frac{d\theta}{\theta}.
\end{equation}
The result now emerges from the well-known, elementary Hardy inequalities (see \cite{hardy},
theorems 327, 328).
\end{proof}
The next result establishes the relevance of $T$ for the Nyman-Beurling approach.
\begin{proposition}\label{rangeofT}
For any $p\in (1,\infty)$, and an interval $E\subseteq (0,\infty)$ the range of $T$ satisfies
\begin{equation}\label{TLpE}
\overline{TL_p(E)}^{L_p}=\overline{\A_E}^{L_p}.
\end{equation}
\end{proposition}
\begin{remark}
For every $f\in L_p(0,\infty)$, $Tf$ is continuous, so the closure operation on the left-hand side above is
necessary. However, for the purpose immediately at hand of proving Theorem
\ref{maintheoremversion1} we only need
\begin{equation}\label{usefulformaintheorem}
TL_p((0,1])\subset\overline{\B}^{L_p}.
\end{equation}
\end{remark}
\begin{remark}
If $f\in L_p(0,1)$ for some $p\in(1,\infty)$ and $$\int_0^1 f(\theta)d\theta=0,$$ then
$Tf\in\overline{\C}^{L_p}$. This is the case for $f=M_1$ by (iii) in Corollary \ref{furtherongamma}. 
\end{remark}
\begin{proof}[Proof of the Proposition]
Fix $p\in(1,\infty)$. For any bounded interval $[a,b]\subseteq E$
\begin{equation}
T\chi_{[a,b]}(x)=\int_a^b \rho\left(\frac{\theta}{x}\right)\frac{d\theta}{\theta}, 
\end{equation}
is a proper Riemann integral for each $x>0$. Let $\theta_{n,k}:=a+(b-a)(k/n)$, and 
\begin{equation}
s_n(x):=\frac{b-a}{n}\sum_{k=1}^{n} \frac{1}{\theta_{n,k}}\rho\left(\frac{\theta_{n,k}}{x}\right).
\end{equation}
The Riemann sums $s_n(x)\in\A$ and $s_n(x)\rightarrow T\chi_{[a,b]}(x)$for each $x>0$. Furthermore it
is trivial to see that $s_n(x)\leq (b-a)/a$ for all $x>0$, whereas $s_n(x)=(b-a)/x$
when $x>b$, so that 
$$s_n(x)\leq \frac{b-a}{a}\chi_{(0,b]}(x)+\frac{b-a}{x}\chi_{(b,\infty)}(x).$$
Hence $\|s_n-T\chi_{[a,b]}\|_p\rightarrow0$. By Proposition \ref{Tispp} we conclude
\begin{equation}\label{Tchiabokay}
T\chi_{[a,b]}\in\overline{\A_E}^{L_p},
\end{equation}
which the time honored density argument and the continuity of $T$ convert into
(\ref{usefulformaintheorem}), and, a fortiori, $\overline{TL_p(E)}^{L_p}\subseteq\overline{\A_E}^{L_p}$. To
finish the proof of (\ref{rangeofT}) we need to show that each function $\rho(\alpha/x)$, $\alpha\in E$ is in
$\overline{TL_p(E)}^{L_p}$. This is achieved as follows: For $\alpha\not=a$ take $1>h\downarrow0$.
Clearly
\begin{equation}\label{differentiateTchi}
\frac{1}{\max(\alpha,x)}\geq \frac{1}{\alpha h}\int_{\alpha(1-h)}^{\alpha}
\rho\left(\frac{\theta}{x}\right)\frac{d\theta}{\theta}\rightarrow
\rho\left(\frac{\alpha}{x}\right), \ \ (\textit{a.e.\ x}).
\end{equation}
By (\ref{Tchiabokay}) $(\alpha h)^{-1}T\chi_{[\alpha(1-h),\alpha]}\in\overline{\A_E}^{L_p}$, and the
above inequalities show it converges in $L_p$-norm to the function $\rho(\alpha/x)$.
If $\alpha=a$ the modification to the above proof is obvious.
\end{proof}

We next introduce the essential, elementary identity.
\begin{lemma}\label{mobiusidentity1theorem}
For every $x>0$
\begin{equation}\label{mobiusidentity1}
\chi_{(1,\infty)}(x)\log x = \int_1^x M(t)\left[\frac{x}{t}\right]\frac{dt}{t}.
\end{equation}
\end{lemma}
\begin{proof}
Here we denote $\chi(S)=1$ if the statement $S$ is true, otherwise $\chi(S)=0$. We start from the well-known
elementary identity
\begin{equation}\label{sumM(n/k)is1}
\chi_{[1,\infty)}(t)=\sum_{n=1}^{\infty}M\left(\frac{t}{n}\right),
\end{equation}
which we multiply by $1/t$ and integrate thus:
\begin{eqnarray}\nonumber
\chi_{[1,\infty)}(x)\log x&=&\int_{0}^{x}\sum_{n=1}^{\infty}M\left(\frac{t}{n}\right)\frac{dt}{t}\\\nonumber
&=&\sum_{n=1}^{\infty}\int_{0}^{x}M\left(\frac{t}{n}\right)\frac{dt}{t}\\\nonumber
&=&\sum_{n=1}^{\infty}\int_{0}^{\frac{x}{n}}M(t)\frac{dt}{t}\\\nonumber
&=&\int_{0}^{x}M(t)\sum_{n=1}^{\infty}\chi\left(t\leq\frac{x}{n}\right)\frac{dt}{t}\\\nonumber
&=&\int_{0}^{x}M(t)\left[\frac{x}{t}\right]\frac{dt}{t}.
\end{eqnarray}
\end{proof}
\begin{proposition}\label{mobiusidentity2theorem}
For every $x>0$ the following identity holds true as an absolutely convergent integral,
without any assumptions on the $L_p$-norms of $M_1$.
\begin{equation}\label{mobiusidentity2}
\lambda(x)=
\int_0^1 M_1(\theta)\rho\left(\frac{\theta}{x}\right)\frac{d\theta}{\theta}.
\end{equation}
\end{proposition}
\begin{proof}
The upper limit of integration in (\ref{mobiusidentity1}) can trivially be substituted by $\infty$, so
we get
\begin{eqnarray}\nonumber
\chi_{(1,\infty)}(x)\log x &=& \int_1^{\infty}
M(t)\left(\frac{x}{t}-\rho\left(\frac{x}{t}\right)\right)\frac{dt}{t}\\\nonumber
 &=&-\int_1^{\infty} M(t)\rho\left(\frac{x}{t}\right)\frac{dt}{t}
\end{eqnarray}
from (iii) in Corollary \ref{furtherongamma} and the (absolute) convergence of the last integral, again due
to $M(t)\ll t(\log t)^{-2}$. Now make the change of variables $t=1/\theta$, and in the formula
obtained substitute $x\mapsto1/x$.
\end{proof}
\begin{remark}\label{mobiusidentity2andtnp}
In \cite{PNTequivlbd} we show that the existence of 
$$
\lim_{\epsilon\rightarrow 0}
\int_{\epsilon}^{1}
M_1(\theta)\rho\left(\frac{\theta}{x}\right)\frac{d\theta}{\theta},
$$
is elementarily equivalent to the prime number theorem.
\end{remark}

\noindent We can now state and prove the main theorems of this paper.

\begin{theorem}[Arithmetical Nyman-Beurling Theorem, I]\label{maintheoremversion1}
The following statements are true for all $p\in(1,\infty)$.\\
\ \\
(a) $\|M_1\|_p<\infty$ implies $\lambda\in\overline{\B}^{L_p}$.\\
\ \\
(b) $\lambda\in\overline{\B}^{L_p}$ implies $\|M_1\|_r<\infty$ for all $r\in(1,p)$.  
\end{theorem}

\begin{proof}[Proof of (a)]
If $\|M_1\|_p<\infty$, then 
$\lambda=TM_1\in TL_p(0,1)\subset\overline{\A_{(0,1)}}^{L_p}=\overline{\B}^{L_p}$ by
(\ref{mobiusidentity2}), Lemma \ref{Tispp}, and (\ref{usefulformaintheorem}). 
\end{proof}
\begin{proof}[Proof of (b)]
If $\lambda\in\overline{\B}^{L_p}$, then $\chi\in\overline{\B}^{L_p}$ as remarked in
(\ref{lambdaimplieschi}). Then by the easy sufficiency part of the Nyman-Beurling Theorem
\ref{extendednymanbeurlingtheorem} $\zeta(s)\not=0$ for $\Re s>1/p$, and this implies by Proposition
\ref{arithmeticalRHp} that $\|M_1\|_r<\infty$ for all $r\in(1,p)$.  
\end{proof}
\begin{remark}
The proof of (a) is elementary, and, interestingly, it corresponds to the ``hard"
\textit{necessity part} of the Nyman-Beurling Theorem \ref{extendednymanbeurlingtheorem}. Note
however that the \textit{strong} form of (a) is connected with the fact that the hypothesis implies by
(\ref{normpMfinite}) that there are no zeroes of $\zeta(s)$ in $\Re s\geq 1/p$.
On the other hand (b), a \textit{weak} statement, corresponding to the ``easy" \textit{sufficiency
part} of the Nyman Beurling Theorem \ref{extendednymanbeurlingtheorem}, is proved essentially by the
traditional argument. 
\end{remark}

The second version is predicated on the new natural approximation $G_n$ defined in (\ref{Gn}). It is
easy to see from Theorem \ref{rangeofT} that $G_n\in\overline{\A_{(1/n,1]}}^p$ for all
$p\in(1,\infty)$.

\begin{theorem}[Arithmetical Nyman-Beurling Theorem, II]\label{maintheoremversion2}
The following statements are true for all $p\in(1,\infty)$.\\
\ \\
(c) $\zeta(s)\not=0$, $\Re s>1/p$ implies $\|G_n-\lambda\|_r\rightarrow0$ for all $r\in(1,p)$.\\
\ \\
(d)  $\|G_n-\lambda\|_p\rightarrow0$ implies $\zeta(s)\not=0$, $\Re s\geq 1/p$.      
\end{theorem}
\begin{proof}[Proof of (c)]
If $\zeta(s)\not=0$ for $\Re s>1/p$, then, by Proposition \ref{arithmeticalRHp}, $\|M_1\|_r<\infty$ for
all $r\in(1,p)$. It is then clear that $G_n=T(M_1\chi_{(1/n,1]})\overset{L_r}{\rightarrow}\lambda$
by the $L_p$-continuity of $T$ (Lemma \ref{Tispp}). 
\end{proof}
\begin{proof}[Proof of (d)]
We proceed by contradiction. Assume there is $s_0$ such that $\zeta(s_0)=0$ and $\Re s_0=1/p_1\geq
1/p$. Therefore $\gamma(n)\not=o(n^{-1/q_1})$ by Corollary \ref{arithmeticalfunctionbehavior}. Now by the
definition (\ref{Gn}) of $G_n$ and (ii) in Lemma \ref{furtherongamma} we have
\begin{eqnarray}\label{normGnminuslambda}
\|G_n-\lambda\|_p^p &=& \int_0^{\infty}\left|\int_0^{1/n}M_1(\theta)\rho\left(\frac{\theta}{x}\right)
\frac{d\theta}{\theta}\right|^p dx\\\nonumber
 &\geq&\int_{1/n}^{\infty}\left|\int_0^{1/n}M_1(\theta)\rho\left(\frac{\theta}{x}\right)
\frac{d\theta}{\theta}\right|^p dx\\\nonumber
 &=& \int_{1/n}^{\infty}\left|\frac{1}{x}\int_0^{1/n}M_1(\theta) d\theta\right|^p dx\\\nonumber
 &=& (p-1)^{-1}|\gamma(n)|^p n^{p-1}, 
\end{eqnarray}
so that
\begin{equation}\label{Gncontradiction}
\|G_n-\lambda\|_p\geq (p-1)^{-1/p}|\gamma(n)|n^{1/q_1}\not\rightarrow0.
\end{equation}
\end{proof}

\begin{remark}
The proof of (c), a weak statement corresponding to the``hard" \textit{necessity part} of the
Nyman-Beurling Theorem \ref{extendednymanbeurlingtheorem}, is easy and quasi elementary.
On the other hand the proof of (d), a strong statement, corresponding to the ``easy" \textit{sufficiency part} of
the Nyman Beurling Theorem, is rather easy, but not elementary.
\end{remark}

\begin{remark}
At least formally one can apply the operators
$$\mathcal D_h: = \frac{1}{h}(K_{(1+h)}-I)$$
to (c), and let $h\downarrow0$ to obtain the corresponding Balazard-Saias result for $S_n$ in \cite{notes1}.
The difficulty in formalizing this argument stems from the fact that, for $\mathcal D_h$ as an operator 
from $L_p$ to itself,  
$\|\mathcal D_h\|\rightarrow\infty$, except when $p=1$. A rigorous proof would be desirable. 
\end{remark}

$G_n$ also behaves nicely pointwise and in $L_1$, as the original natural approximations. 
To see this we first need a lemma.
\begin{lemma}\label{speciallemma}
For $\theta>n$
\begin{equation}
\int_n^{\infty}\rho\left(\frac{x}{\theta}\right)\frac{dx}{x^2}\ll \frac{\log \theta}{\theta}.
\end{equation}
\end{lemma}
\begin{proof}[Proof of the Lemma]
For $\theta>n$ we have
\begin{eqnarray}\nonumber
\int_n^{\infty}\rho\left(\frac{x}{\theta}\right)\frac{dx}{x^2}&=&
\frac{1}{\theta}\int_n^{\theta}\frac{dx}{x}+
\int_{\theta}^{\infty}\rho\left(\frac{x}{\theta}\right)\frac{dx}{x^2}\\\nonumber
 &\leq& \frac{\log \theta}{\theta}+\frac{1}{\theta}.
\end{eqnarray}
\end{proof}

\begin{proposition}\label{Gnisnatural}
$G_n$ satisfies these properties
\begin{eqnarray}
G_n(x)-\lambda(x) &\rightarrow& 0, \ \ (\forall x>0),\\\label{Gnconvergespointwise}
 \int_0^1 |G_n(x)-\lambda(x)|dx &\rightarrow& 0.\label{GnconvergesinL1}
\end{eqnarray}
\end{proposition}
\begin{proof}
The first statement follows easily from the fact that the integral in (\ref{mobiusidentity2}) is
absolutely convergent. Changing variables in the first iterated integral below we get
\begin{equation}
\|G_n-\lambda\|_1=\int_0^1\left|\int_0^{1/n}M_1(\theta)\rho\left(\frac{\theta}{x}\right)
\frac{d\theta}{\theta}\right|dx
=\int_1^{\infty}\left|\int_n^{\infty}M(\theta)\rho\left(\frac{x}{\theta}\right)
\frac{d\theta}{\theta}\right|\frac{dx}{x^2}.
\end{equation}
Now we split the outer integral on the righthand side as $\int_1^n + \int_n^{\infty}$. The first one
easily evaluates to $|\gamma(n)|\log n$ taking into account (ii) and (iii) in Lemma
\ref{furtherongamma}. This term converges to zero in view of (\ref{Mggammanatural}) and an elementary
error term in the prime number theorem. The second one is bounded by
\begin{eqnarray}\nonumber
\int_n^{\infty}\int_n^{\infty}|M(\theta)|\rho\left(\frac{x}{\theta}\right)
\frac{d\theta}{\theta}\frac{dx}{x^2} &=&\int_n^{\infty}\frac{|M(\theta)|}{\theta}
\left(\int_n^{\infty}\rho\left(\frac{x}{\theta}\right)\frac{dx}{x^2}\right)d\theta\\\nonumber
 &\ll & \int_n^{\infty}|M(\theta)|\frac{\log \theta}{\theta^2}d\theta\\\nonumber
 &\ll & \int_n^{\infty}\frac{d\theta}{\theta\log^2 \theta}\rightarrow0,\ \ (n\rightarrow\infty),
\end{eqnarray}
where we have applied in succession Fubini's theorem, Lemma \ref{speciallemma}, and an elementary
error tem for the prime number theorem of the form $M(x)\ll x(\log x)^{-3}$. 
\end{proof}

\section{On divergence of certain natural approximations}\label{divergenceofnaturalapproximations}
Throughout this section $1<p<\infty$. All natural approximations considered converge both a.e. and in $L_1$
either to $\lambda$ or to $-\chi$, hence not converging in $L_p$ to the corresponding generator is equivalent
to diverging in $L_p$.

\subsection{Divergence of approximations to $\lambda$.} The main result needs no proof as it is just the
counterpositive of statement (d) in Theorem \ref{maintheoremversion2}, namely:

\begin{proposition}\label{Gndivergence}
If $\zeta(s)$ has a zero with real part $\geq 1/p$, then $G_n$ diverges in $L_p$. In particular $G_n$
diverges in $L_2$.
\end{proposition}
\begin{remark}\label{analyzeweaktostrong}
This proposition shows that \textit{in general} the \textit{weak} implication,\\
\ \\
(c) $\zeta(s)\not=0$, $\Re s>1/p$ implies $\|G_n-\lambda\|_r\rightarrow0$ for all $r\in(1,p)$\\
\ \\ 
in Theorem \ref{maintheoremversion2} \textit{cannot be made stronger} to include $r=p$. The hypothesis of
(c) can hold only if $1\leq p \leq 2$. Although we resolved at the outset to keep $1<p<\infty$, our
resolve is weak, so we note that for $p=1$ the strong version is true because of Theorem
\ref{Gnisnatural}. For $p=2$ the strong statement is definitely false for there are zeroes on $\Re s=1/2$.
In the case $1<p<2$ a simple logical analysis shows that the only interesting case is $\beta = 1/p$. Now,
either there are roots on the line $\Re s = \beta$, then the strong statement is false; or else, there are no roots
on that line, then we can say nothing at present. This is related to Question \ref{betaconjecture}.
\end{remark}

\begin{remark}
By Corollary \ref{arithmeticalfunctionbehavior} there is a subsequence of \textit{zero-crossings} of $\gamma(n)$ where
clearly $|\gamma(n)|<1/n$. For this subsequence the contradiction (\ref{Gncontradiction})
would not hold. Thus the possibility remains open that there is a subsequence of $G_n$ that does converge.
This peculiarity is common to all natural approximations discussed here. But there are reasons to believe
this is a mirage.   
\end{remark}

To probe a little into the possible mirage we now bring to bear the existence of an isometry\footnote{It
is actually a unitary operator, but that is not relevant here.} of $L_2(0,\infty)$ denoted by U in \cite{IUO}
satisfying the following conditions, where we let $\rho_1(x)=\rho(1/x)$:

\begin{eqnarray}\label{Uinvariant}
U K_a &=& K_a U,\ (a>0),\\\label{Urho1}
U\rho_1(x) &=& \frac{\rho(x)}{x},\\\label{Uchi}
U\chi(x) &=& \frac{\sin(2\pi x)}{\pi x}.
\end{eqnarray}
For $f\in \A$ as in (\ref{abcfunction})
\begin{equation}\label{Uf}
Uf(x)=\frac{1}{x}\sum_{k=1}^n c_k \theta_k \rho\left(\frac{x}{\theta_k}\right).
\end{equation}
If we apply this to the Riemann sums of $Tf$, when $f$ is continuous of compact support, and make the obvious
modifications to the reasoning in Lemma \ref{Tispp} and Proposition \ref{rangeofT}, we obtain:

\begin{lemma}\label{UT} 
For $f\in L_2(0,\infty)$
\begin{equation}\label{UTf}
UTf(x)=\frac{1}{x}\int_0^{\infty} f(\theta) \rho\left(\frac{x}{\theta}\right)d\theta.
\end{equation}
Moreover the right-hand side defines a continuous extension to all $L_p(0,\infty)$.
\end{lemma}

\begin{remark}\label{UnotLpbutUTLp}
At present we shall use this lemma only in $L_2$. It is however interesting to see how $UT$ extends to all $L_p$'s
given the fact that $U$ cannot be extended continuously to any $L_p$ other than for $p=2$ (see \cite{Utypeset}).
When restricted to $f\in L_2(0,1)$ the integral of the right-hand side of (\ref{UTf}) is the
Hilbert-Schmidt operator studied by J. Alc\'{a}ntara-Bode in \cite{alcantara} where it is shown at the outset
that the Riemann Hypothesis is equivalent to the injectivity of this operator.
\end{remark}
\noindent The above lemma leads to the simple calculation:
\begin{equation}\label{UGnnearzero}
UG_n(x)=H_2(n),\ \textnormal{for } x\in(0,1/n),
\end{equation}
which spells further trouble for the $L_2$ convergence of subsequences of $G_n$:
\begin{proposition}\label{Gnmoreproblems}
\begin{equation}\label{Gndoublejeopardy}
\|\lambda-G_n\|_2\gg \max\left(n^{1/2}\gamma(n),n^{1/2}H_2(n)\right).
\end{equation}
\end{proposition}
\begin{remark}
Since there are roots of $\zeta(s)$ on $\Re s=1/2$, neither $n^{1/2}\gamma(n)$ nor $n^{1/2}H_2(n)$ converge
to zero by Corollary \ref{arithmeticalfunctionbehavior}, and most likely they are unbounded as
$n\rightarrow\infty$. However, optimism about almost periodicity of these functions may induce the idea that
their zero crossings implied also by Corollary \ref{arithmeticalfunctionbehavior} will be close together an
infinite number of times.  
\end{remark}
\begin{proof}[Proof of Proposition \ref{Gnmoreproblems}]
That $\|\lambda-G_n\|_2\gg n^{1/2}\gamma(n)$ is simply (\ref{normGnminuslambda}) for $p=2$. For the second part we use
(\ref{UGnnearzero}):
\begin{eqnarray}\label{UGnUlambda}
\|G_n-\lambda\|_2^2 &=& \|UG_n-U\lambda\|_2^2\\\nonumber
                    &\geq& \int_0^{1/n}\left|UG_n(x)-U\lambda(x)\right|^2 dx\\\nonumber
                    &\geq& \int_0^{1/n}\left|H_2(n)-U\lambda(x)\right|^2 dx\\\nonumber
                    &\gg& n^{-1}|H_2(n)|.
\end{eqnarray}
\end{proof}
A finer analysis of selected intervals in $(1/n,\infty)$ seems likely to produce an infinite number of
barriers increasing the lower bound in (\ref{UGnUlambda}), so that one may be inclined to think that
all subsequences of $G_n$ diverge in $L_2$.\\
\ \\
An even more natural looking approximation of $\lambda$ is obtained by writing the simplest
Riemann sum of the integral (\ref{mobiusidentity2}), namely
$$R_n(x):=\sum_{k=1}^{n-1} \frac{1}{k}M\left(\frac{n}{k}\right)\rho\left(\frac{k}{nx}\right),$$
which happens to be a Beurling function in $\C$ with an uncanny resemblance to the dual approximation
$F_n$ defined by (\ref{Fn}). But bear in mind that the integral (\ref{mobiusidentity2}) is not a proper Riemann integral,
and we have not yet been able to show that $R_n$ is a natural approximation, in the sense that it satisfies a
weak Beurling theorem such as Theorem \ref{maintheoremversion2}, so we state the following true theorem without
proof: 

\begin{proposition}\label{Rndivergence}
$R_n$ diverges in $L_2$
\end{proposition}

Yet another approximation could be defined by truncation, say

$$T\left(\min\left(n,\max(M_1,-n)\right)\right).$$
\ \\
We shall not pursue this matter here either, but it seems to deserve some attention.

\subsection{Divergence of approximations to $-\chi$}
We may treat $S_n$ and $V_n$ together, defined in (\ref{Sn}), (\ref{Vn}), since $\|S_n-V_n\|_2\rightarrow0$ 
Here is then the corresponding divergence result for $S_n$.

\begin{proposition}\label{Sndivergence}
If there is some zero of $\zeta(s)$ with real part $1/p$ then $S_n$
and $V_n$ diverge in $L_p$. In particular $S_n$ and $V_n$ diverge in $L_2$.
\end{proposition}
\begin{proof}
The hypothesis on the zero of $\zeta(s)$ implies, by Corollary \ref{arithmeticalfunctionbehavior}, that 

\begin{equation}\label{gorder}
g(n)\not=o(n^{-1/q}).
\end{equation}
Now assume by contradiction that
$S_n$ converges in
$L_p$, so it must converge to $-\chi$. On the other hand, noting that $kx>1$ when $x>1/m$ and $k>m$, we get

\begin{eqnarray}\nonumber
\|S_n-S_m\|_p^p &\geq& 
\int_{1/m}^{\infty}\left|\sum_{k=m+1}^n \mu(k)\rho\left(\frac{1}{kx}\right)\right|^p dx\\\nonumber
 &=& \frac{1}{p-1}m^{p-1}|g(n)-g(m)|^p.
\end{eqnarray}
Then letting $n\rightarrow\infty$ we obtain

\begin{equation}\label{hurdle1}
\|\chi+S_m\|_p^p \geq \frac{1}{p-1}m^{p-1}|g(m)|^p.
\end{equation}
Since the left-hand side goes to zero when $m\rightarrow\infty$ this contradicts (\ref{gorder}).
\end{proof}

\begin{remark}
The above proposition implies that \textit{in general} the weak implication of Balazard-Saias ((i) implies
(vii) in \cite{notes1}, see also \cite{lee})\\
\ \\ 
$\zeta(s)\not=0$, $\Re s>1/p$ implies $\|S_n+\chi\|_r\rightarrow0$ for all $r\in(1,p)$\\
\ \\
cannot be made stronger to include $r=p$. An analysis analogous to that carried out
for $G_n$ in Remark \ref{analyzeweaktostrong} is possible here too. Mutatis mutandis the conclusions are the
same. But a cautionary note is in order. We have not been able to treat the $L_p$ case for $B_n$, other than
for $p = 1$ or $2$.
\end{remark}

\begin{remark}
Again, the existence of a subsequence of zero-crossings of $g(n)$ given by Corollary
\ref{arithmeticalfunctionbehavior} indicates that this subsequence is still a candidate in the running to converge in
$L_p$-norm to $-\chi$. However as with $G_n$ we now prove a stricter failure for $S_n$ in the $L_2$ case. 
\end{remark}

\begin{proposition}
There exists a constant $C>0$ such that
\begin{equation}\label{hurdle1and2}
\|\chi+S_{n}\|_{2}\geq \max\left(\frac{C}{\sqrt{n}}|M(n)+2|,|g(n)|\sqrt{n}\right).
\end{equation}
\end{proposition}
\begin{proof}
For $p=2$ inequality (\ref{hurdle1}) translates into

\begin{equation}\label{hurdle1norm2}
\|\chi+S_n\|_2 \geq |g(n)|\sqrt{n} .
\end{equation}
On the other hand if we apply $U$ to $S_n$ we get
\begin{equation}
US_n(x)=M(n),\ \ (0<x<1/n).
\end{equation}
Hence

\begin{equation}\label{normofUchiplusUSn}
\|\chi+S_{n}\|_{2}^2\geq\int_0^{1/n}\left|\frac{1}{\pi x}\sin(2\pi x)+M(n)\right|^2 dx,
\end{equation}
and

\begin{equation}\label{SndivergesinL2}
\|\chi+S_{n}\|_{2}\geq C \frac{1}{\sqrt{n}}|M(n)+2|,
\end{equation}
for some positive constant $C$.
\end{proof}
\noindent Odlyzko and te Riele \cite{odlyzko} have conjectured that

$$\limsup_{n\rightarrow\infty} \frac{|M(n)|}{\sqrt{n}} = \infty,$$
in which case $S_n$ would not even be bounded in $L_2$, endangering also the possibility of
a strong version of condition (vi) in Balazard-Saias's work \cite{notes1}. On the other hand, by
Corollary \ref{arithmeticalfunctionbehavior} there is a subsequence where $M(n)=-2$, and we know there
is a subsequence where $g(n)$ crosses zero, with $g(n)\leq1/n$. Nevertheless, as for $G_n$, one may suspect that
there is no $L_2$-convergent subsequence of $S_n$.

The initial natural approximation $B_n$ is more resilient. We already remarked that it is not equivalent to
$S_n$, neither is it a \textit{series} as defined below. The fact that $B_n(x)=-1$ in $(1/n,1)$ destroys the
possibility of using the same argument of Proposition \ref{Sndivergence}. However with the help of the
operator
$U$ we can dispose of the $L_2$-case both for $B_n$ and $F_n$.

\begin{proposition}\label{BnFndivergeinL2}
Neither $B_n$ nor $F_n$ converge in $L_2$.
\end{proposition}
\begin{proof}
The $U$ defining properties (\ref{Uinvariant}), (\ref{Urho1}), as well as (\ref{Mggammanatural}) give
$B_n(x)=-n\gamma(n)$ in $(0,1/n)$. Assume by contradiction that $B_n$ converges in $L_2$, then so does
$UB_n$, and therefore

$$0\leftarrow\int_0^{1/n}|UB_n(x)|^2 dx = n|\gamma(n)|^2,$$
which contradicts Corollary \ref{arithmeticalfunctionbehavior}. Likewise the wholly analogous computation
$UF_n(x)=M(n)-1$ in $(0,1/n)$ yields the divergence of $F_n$ in $L_2$.
\end{proof}

\noindent Analogous considerations as for $S_n$ apply in relation to the possibility of a subsequence of
$B_n$ or of $F_n$ converging in $L_2$.\\
\ \\

To round off the presumption of divergence of the natural approximations in $L_2$, we prove
Proposition \ref{noseriesconverges}, a result, suggested by M. Balazard\footnote{Personal communication.},
stating that no \textit{series} of a certain kind in $\C$ can converge to $-\chi$ in $L_2$.
 
Denote by $\C^{nat}$
the subspace generated by the linearly independent functions $\{e_{k}| k \geq 2\}$, where
\begin{equation}
e_{k}(x):=\rho\left(\frac{1}{kx}\right)-\frac{1}{k}\rho\left(\frac{1}{x}\right).
\end{equation}
Note that

\begin{equation}
V_{n}=\sum_{k=2}^{n}\mu(k)e_{k}.
\end{equation}
A \textit{series} in $\C^{nat}$ is defined as any sequence of type

\begin{equation}\label{seriesinC}
f_{n}=\sum_{k=2}^{n}c_{k}e_{k},\ \ \ (n\geq2).
\end{equation}
We can now state:

\begin{proposition}\label{noseriesconverges}
No series in $\C^{nat}$ converges in $L_{p}(0,1)$ to $-\chi$ if there is a zero of $\zeta(s)$ with
real part $1/p$. In particular, no series in $\C^{nat}$ converges to $-\chi$ in $L_2(0,1)$. 
\end{proposition}

\noindent To achieve the proof of this theorem we need a lemma.

\begin{lemma}\label{coeffsconverge}
Let $f_{n}$ be a sequence in $\C^{nat}$ converging pointwise to $-\chi$. Assume $f_{n}$ is written as
\begin{equation}
f_{n}(x)=\sum_{k=1}^{n}a_{n,k}\hspace{.5mm}\rho\left(\frac{1}{kx}\right),
\end{equation}
then
\begin{equation}\label{coeffkconvergestomuk}
a_{n,j}\rightarrow\mu(j),\ \ (n\rightarrow\infty),
\end{equation}
for every $j \geq 1$.
\end{lemma}

\begin{proof}
Each $f_n\in\C$, so condition (\ref{sumckthetakiszero}) implies it is the right-continuous, step function
\begin{equation}
f_n(x)=-\sum_{k=1}^n a_{n,k}\left[\frac{1}{kx}\right],
\end{equation}
which is constant on every interval $$\left(\frac{1}{j+1},\frac{1}{j}\right],\ \ j=1, 2,
\dots.$$
Therefore pointwise convergence trivially implies
\begin{equation}\label{system}
-\lim_{n\rightarrow\infty}f_n(1/j)=\lim_{n\rightarrow\infty}\sum_{k=1}^n a_{n,k}\left[\frac{j}{k}\right]
\rightarrow1,\ \ j=1, 2,
\dots.
\end{equation}
Now we proceed by induction. For $j=1$ it is clear that (\ref{system}) gives $a_{n,1}\rightarrow1=\mu(1)$.
Next assume for $j>1$ that $a_{n,k}\rightarrow\mu(k)$ for $1\leq k \leq j-1$, then the limit in
(\ref{system}) yields
$$\sum_{k=1}^{j-1}\mu(j)\left[\frac{j}{k}\right] + a_{n,j}\rightarrow 1.$$
But comparing this to the well-known
$$\sum_{k=1}^{j}\mu(j)\left[\frac{j}{k}\right] = 1,$$
we obtain the desired $a_{n,j}\rightarrow\mu(j)$ as $n\rightarrow\infty$.
\end{proof}

\begin{remark}
In some sense this Lemma shows the \textit{inevitability} of the natural approximation $S_n$.  
\end{remark}

\begin{proof}[Proof of Proposition \ref{noseriesconverges}]
We have trivially
\begin{equation}
f_{n}(x)=
-\left(\sum_{k=2}^{n}\frac{c_{k}}{k}\right)\rho\left(\frac{1}{x}\right)
+\sum_{k=2}^{n}c_{k}\hspace{.5mm}\rho\left(\frac{1}{kx}\right).
\end{equation}
Assume by contradiction that $\|\chi+f_{n}\|_{p}\rightarrow0$. For the step functions involved this clearly implies
pointwise convergence, then, from Lemma \ref{coeffsconverge}, we get
$c_{k}=\mu(k)$ for each $k\geq 2$, which, by the way, forces (\ref{coeffkconvergestomuk}) to hold
for $k=1$ too. But this immediately implies that $f_{n}=V_{n}$. However, $V_n$ diverges in $L_p$ by Proposition
\ref{Sndivergence}, so we have obtained a contradiction.
\end{proof}
\ \\
\ \\
\textbf{Ackowledgements.} We thank M. Balazard and E. Saias for useful conversations, and A. M. Odlyzko for his
generous help in the matters treated in Section \ref{lemmae}.

\bibliographystyle{amsplain}

\ \\
\noindent Luis B\'{A}EZ-DUARTE\\
Departamento de Matem\'{a}ticas\\
Instituto Venezolano de Investigaciones Cient\'{\i}ficas\\
Apartado 21827, Caracas 1020-A\\
Venezuela\\

\end{document}